\documentclass[12pt]{amsart}
\usepackage{amssymb}
\usepackage{cite}
\usepackage[textwidth=13.5cm,textheight=21cm]{geometry}

\newtheorem{theorem}{Theorem}[section]
\newtheorem{proposition}[theorem]{Proposition}
\newtheorem{lemma}[theorem]{Lemma}
\newtheorem{corollary}[theorem]{Corollary}
\theoremstyle{remark}
\newtheorem{remark}[theorem]{Remark}

\newcommand{\Li}{\operatorname{Li}}
\newcommand{\li}{\operatorname{li}}

\title[The centre of mass of the primes]{The centre of mass of the primes and a reflection of Chebyshev's bias}

\author{Elham Mehdi-Nezhad}
\address{Department of Mathematics and Applied Mathematics, University of the Western Cape, Private Bag X17, Bellville 7535, Cape Town, South Africa}
\email{emehdinezhad@uwc.ac.za}

\author{Amir M. Rahimi}
\address{School of Mathematics, Institute for Research in Fundamental Sciences (IPM), P.O. Box 19395-5746, Tehran, Iran}
\email{amrahimi@ipm.ir}


\numberwithin{equation}{section}
\begin{document}

\begin{abstract}
For $x\ge2$ we place a unit mass at $\log p$ for each prime $p\le x$, and let $D(x)$ be the distance from the centre of mass of this system to $\log x$; thus $D(p_n)=\log(p_n/G_n)$, where $G_n$ is the geometric mean of the first $n$ primes. Let $D_0(x)$ be the distance obtained when $\pi$ is replaced by $\operatorname{Li}(x)=\int_2^x dt/\log t$. We prove that $D(x)-D_0(x)\ll x^{\Theta-1+\varepsilon}$, where $\Theta$ is the supremum of the real parts of the zeros of $\zeta(s)$, and that the exponent cannot be lowered; in particular, the Riemann hypothesis is equivalent to $\log(p_n/G_n)=D_0(p_n)+O_\varepsilon(p_n^{-1/2+\varepsilon})$ for every $\varepsilon>0$. Under the Riemann hypothesis we prove a pointwise mirror theorem: for large $x$ the centre of mass lies to the left of its predicted position only when $\pi(x)<\li(x)-1.95\sqrt x/\log x$. Under the Riemann hypothesis and the linear independence hypothesis, the set of such $x$ has the same logarithmic density as the set on which $\pi(x)>\li(x)$, although for large $x$ the two sets are disjoint. We prove analogues for primes in arithmetic progressions, for the higher moments of the system, and, under a linear independence hypothesis, for $\mathbb F_q[T]$, where the biases of even and odd degrees are interchanged.
\end{abstract}

\maketitle

\section{Introduction}

The geometric mean $G_n=(p_1\cdots p_n)^{1/n}$ of the first $n$ primes satisfies $p_n/G_n\to e$. This follows from the prime number theorem; it was conjectured by Vrba and proved by S\'andor and Verroken \cite{SV}. S\'andor \cite{S12} showed that $p_n/G_n>e$ for $n\ge10$, and Hassani \cite{Hassani13}, Kourbatov \cite{Kourbatov} and Axler \cite{Axler18} obtained further expansions and explicit bounds. These results concern the asymptotic expansion of $p_n/G_n$. Here we study a different question: the oscillation of $\log(p_n/G_n)$ about a model, and its relation to the zeros of $\zeta(s)$ and to Chebyshev's bias.

We regard the primes as a system of point masses: for $x\ge2$ we place a unit mass at $\log p$ for each prime $p\le x$. The total mass is $\pi(x)$, the first moment about the origin is $\theta(x)=\sum_{p\le x}\log p$, and the density of the masses increases towards the point $\log x$. Let
\[
D(x)=\log x-\frac{\theta(x)}{\pi(x)}
\]
be the distance from the centre of mass to $\log x$, so that $D(p_n)=\log(p_n/G_n)$. The moment of the system about $\log x$ is the Riesz mean $F(x)=\int_2^x\pi(t)\,dt/t$ (Lemma~\ref{lem:identity}), so that $D=F/\pi$ combines a smoothed and an unsmoothed counting function. Seen from $\log x$, the masses are distributed approximately according to the exponential law, so that $D(x)\to1$, which is the statement $p_n/G_n\to e$.

Let $\Li(x)=\int_2^x dt/\log t$, which differs from the logarithmic integral $\li(x)$ by the constant $\li(2)$. Replacing $\pi$ by $\Li$ in an exact formula for $D$ (Lemma~\ref{lem:identity}) gives the \emph{model distance}
\[
D_0(x)=\frac{1}{\Li(x)}\int_2^x\frac{\Li(t)}{t}\,dt .
\]
Define $\ell_0,\ell_1,\dots$ by the formal identity $\bigl(\sum_{k\ge0}k!\,t^k\bigr)^{-1}=1-\sum_{j\ge0}\ell_jt^{j+1}$, so that $(\ell_j)=(1,1,3,13,71,\dots)$; these are the coefficients of Panaitopol's formula for $\pi(x)$ \cite{Pan00}. Then $D_0(x)=\sum_{j<J}\ell_j(\log x)^{-j}+O_J((\log x)^{-J})$ for every $J\ge1$. Axler \cite[Proposition~2.1]{Axler18} proved that $\log(p_n/G_n)$ has the same expansion, confirming approximations proposed by Kourbatov \cite{Kourbatov}; this also follows at once from \eqref{eq:DD0} below and the prime number theorem with de la Vall\'ee Poussin's error term. The finer question is how well $D_0(p_n)$ itself approximates $\log(p_n/G_n)$.

\begin{theorem}\label{thm:main}
Let $\Theta=\sup\{\operatorname{Re}\rho:\ \zeta(\rho)=0,\ 0<\operatorname{Re}\rho<1\}$.
\begin{enumerate}
\item[(a)] For every $\varepsilon>0$ we have $D(x)-D_0(x)\ll_\varepsilon x^{\Theta-1+\varepsilon}$. If $\Theta>\frac12$, then for every $\beta<\Theta$
\[
\limsup_{x\to\infty}\frac{D(x)-D_0(x)}{x^{\beta-1}}=+\infty,\qquad\liminf_{x\to\infty}\frac{D(x)-D_0(x)}{x^{\beta-1}}=-\infty .
\]
\item[(b)] $D(x)-D_0(x)=\Omega_\pm\bigl(x^{-1/2}\log\log\log x\bigr)$.
\end{enumerate}
In particular, $D(x)-D_0(x)$ changes sign infinitely often, $\inf\{\alpha:\ D(x)-D_0(x)\ll x^{\alpha}\}=\Theta-1$, and the Riemann hypothesis holds if and only if $D(x)=D_0(x)+O_\varepsilon(x^{-1/2+\varepsilon})$ for every $\varepsilon>0$.
\end{theorem}

The last statement is the analogue for the centre of mass of von Koch's theorem \cite{vonKoch} and its classical converse.

\begin{corollary}\label{cor:Gn}
\begin{enumerate}
\item[(a)] The Riemann hypothesis holds if and only if, for every $\varepsilon>0$,
\[
\log\frac{p_n}{G_n}=D_0(p_n)+O_\varepsilon\bigl(p_n^{-1/2+\varepsilon}\bigr).
\]
\item[(b)] $\log(p_n/G_n)-D_0(p_n)=\Omega_\pm\bigl(p_n^{-1/2}\log\log\log p_n\bigr)$; in particular, $\log(p_n/G_n)-D_0(p_n)$ changes sign infinitely often.
\end{enumerate}
\end{corollary}

Unconditionally, \eqref{eq:DD0} and the prime number theorem with the error term of Vinogradov and Korobov give $\log(p_n/G_n)-D_0(p_n)\ll\exp\bigl(-c(\log p_n)^{3/5}(\log\log p_n)^{-1/5}\bigr)$; Corollary~\ref{cor:Gn} shows that a power saving is equivalent to the Riemann hypothesis. For Mertens' product, oscillation results of the type of Theorem~\ref{thm:main}(b) were proved by Diamond and Pintz \cite{DP}. The moment $F$ alone also detects the Riemann hypothesis: it holds if and only if $\int_2^x(\pi(t)-\Li(t))\,dt/t<0$ for all large $x$ (Corollary~\ref{cor:Riesz}).

The sign of $D(x)-D_0(x)$ is governed by Chebyshev's bias, but in reflected form. Let $\gamma_0$ be Euler's constant. Assume the Riemann hypothesis and put
\[
\sigma_0=\sum_\rho\frac{1}{|\rho|^2}=\sum_\rho\frac{1}{\rho(1-\rho)}=2+\gamma_0-\log4\pi=0.0461\ldots,
\]
the sums running over the non-trivial zeros of $\zeta(s)$; the last identity is classical and holds unconditionally (see \cite[Ch.~12]{Davenport}).

\begin{theorem}\label{thm:mirror}
Assume the Riemann hypothesis, and let $\varepsilon>0$. For all sufficiently large $x$,
\[
\begin{aligned}
\pi(x)<\li(x)-(2+\sigma_0+\varepsilon)\frac{\sqrt x}{\log x}\ &\Longrightarrow\ D(x)>D_0(x)\\
&\Longrightarrow\ \pi(x)<\li(x)-(2-\sigma_0-\varepsilon)\frac{\sqrt x}{\log x}.
\end{aligned}
\]
\end{theorem}

Hence, for large $x$, the inequalities $D(x)>D_0(x)$ and $\pi(x)>\li(x)$ never hold together: the centre of mass lies to the left of its predicted position only when $\pi(x)$ is well below $\li(x)$. We call Theorem~\ref{thm:mirror} and its analogue for residue classes (Theorem~\ref{thm:mirrorAP}) mirror theorems. The explanation lies in the explicit formula for $\pi(x)(D(x)-D_0(x))$, in which the zero $\rho=\frac12+i\gamma$ occurs with the weight
\[
\frac1\rho-\frac{1}{\rho^2}=\frac{\rho-1}{\rho^2}=-\frac{\bar\rho}{\rho}\cdot\frac1\rho .
\]
This differs from the weight $1/\rho$ in the explicit formula for $\pi(x)-\li(x)$ by the factor $-\bar\rho/\rho$ of modulus one, while the bias coming from the squares of primes is the same. The constant $2$ in Theorem~\ref{thm:mirror} is the constant term in Lemma~\ref{lem:integral}, which also comes from the squares of primes, and $\sigma_0$ bounds the remaining sum $\sum_\rho x^{i\gamma}/\rho^2$. More precisely, the normalized function $\frac{\log x}{\sqrt x}\pi(x)(D(x)-D_0(x))$ differs from $-2-D_0(x)\frac{\log x}{\sqrt x}(\pi(x)-\li(x))$ by at most $\sigma_0+o(1)$ in absolute value. A set $S\subset[2,\infty)$ has logarithmic density $\delta$ if $\frac1{\log X}\int_{S\cap[2,X]}dx/x\to\delta$ as $X\to\infty$. Let LI denote the hypothesis that the positive ordinates of the zeros of $\zeta(s)$ are linearly independent over $\mathbb Q$. Under LI the numbers $x^{i\gamma}$ behave like independent random variables uniformly distributed on the unit circle, and multiplying each of them by the factor $-\bar\rho/\rho$ does not change their joint law; this explains Theorem~\ref{thm:skewes} below. Rubinstein and Sarnak \cite{RS} showed, under the Riemann hypothesis and LI, that the set of $x$ with $\pi(x)>\li(x)$ has logarithmic density $\delta_1=0.00000026\ldots$.

\begin{theorem}\label{thm:skewes}
Assume the Riemann hypothesis and LI. Then the set of $x$ with $D(x)>D_0(x)$ has logarithmic density $\delta_1$.
\end{theorem}

Theorem~\ref{thm:skewes} thus reflects the symmetry of the limiting distribution of $\frac{\log x}{\sqrt x}(\pi(x)-\li(x))$ about its mean $-1$, and not any relation between the two exceptional sets, which by Theorem~\ref{thm:mirror} are disjoint for large $x$. Theorem~\ref{thm:skewes} is the analogue for the centre of mass of a theorem of Lamzouri \cite{LamM} on Mertens' product: under the Riemann hypothesis and LI, the inequality $\prod_{p\le x}(1-1/p)^{-1}<e^{\gamma_0}\log x$ holds on a set of logarithmic density $\delta_1$; there the zero $\rho$ carries the weight $1/(\rho-1)$, which on the critical line also has modulus $1/|\rho|$.

Similar results hold for primes in arithmetic progressions. For $q\ge3$ and $(a,q)=1$ let $\pi(x;q,a)$ be the number of primes $p\le x$ with $p\equiv a\pmod q$, let $\theta(x;q,a)$ be the sum of $\log p$ over these primes, and for $x$ with $\pi(x;q,a)\ge1$ let
\[
D(x;q,a)=\log x-\frac{\theta(x;q,a)}{\pi(x;q,a)}
\]
be the distance from $\log x$ of the centre of mass of these primes. By GRH we mean the Riemann hypothesis for all Dirichlet $L$-functions modulo $q$, including $\zeta(s)$, and by LI the hypothesis of Rubinstein and Sarnak that the non-negative ordinates of the zeros of the $L(s,\chi)$, $\chi\ne\chi_0$, form a multiset that is linearly independent over $\mathbb Q$. Let $\nu(q,a)$ be the number of solutions of $u^2\equiv a\pmod q$, and for a non-principal character $\chi$ modulo $q$ let $\sigma_\chi=\sum_{\rho_\chi}|\rho_\chi|^{-2}$, the sum running over the non-trivial zeros of the primitive $L$-function that induces $L(s,\chi)$. Under GRH, $\sigma_\chi=2\operatorname{Re}\lambda_1(\chi)$, where $\lambda_1(\chi)=\sum_{\rho_\chi}1/\rho_\chi$ is the first coefficient in Li's criterion \cite{Li,BL}; for $\zeta(s)$ this gives $\sigma_0=2\lambda_1$.

\begin{theorem}\label{thm:mirrorAP}
Let $q\ge3$, assume GRH for the characters modulo $q$, and let $a,b$ be residues prime to $q$ with
\[
\Delta_{a,b}:=2\bigl(\nu(q,a)-\nu(q,b)\bigr)-\sum_{\chi\ne\chi_0}|\chi(a)-\chi(b)|\,\sigma_\chi>0 .
\]
Then for every $\varepsilon>0$ and all sufficiently large $x$, the inequality $D(x;q,a)\ge D(x;q,b)$ implies
\[
\pi(x;q,b)-\pi(x;q,a)>\Bigl(\frac{\Delta_{a,b}}{\varphi(q)}-\varepsilon\Bigr)\frac{\sqrt x}{\log x}.
\]
Consequently, if $\pi(x;q,a)\ge\pi(x;q,b)$, then $D(x;q,b)>D(x;q,a)$.
\end{theorem}

For a primitive character $\chi$ modulo $q^*$ with $\chi(-1)=(-1)^{\mathfrak a}$, $\mathfrak a\in\{0,1\}$, the Riemann hypothesis for $L(s,\chi)$ gives $\sigma_\chi=2\operatorname{Re}\frac{L'}{L}(1,\chi)+\log\frac{q^*}{\pi}+\frac{\Gamma'}{\Gamma}\bigl(\frac{1+\mathfrak a}{2}\bigr)$ (Lemma~\ref{lem:sigma}). For $q=4$ this gives $\sigma_{\chi_{-4}}=\gamma_0+6\log2+5\log\pi-8\log\Gamma(\frac14)=0.1555\ldots$, and with $a=1$, $b=3$ we have $\Delta_{1,3}=4-2\sigma_{\chi_{-4}}=3.688\ldots$. Thus, under GRH, for all large $x$, whenever $\pi(x;4,1)\ge\pi(x;4,3)$ the centre of mass of the primes $\equiv3\pmod4$ is further from $\log x$, and whenever the centre of mass of the primes $\equiv1\pmod4$ is at least as far from $\log x$, we have $\pi(x;4,3)-\pi(x;4,1)>1.84\sqrt x/\log x$. Similarly $\sigma_{\chi_{-3}}=\gamma_0+8\log2-2\log3+7\log\pi-12\log\Gamma(\frac13)=0.1132\ldots$, and $\Delta_{1,2}=3.773\ldots$ for $q=3$ and $\Delta_{1,5}=3.773\ldots$ for $q=6$. On the other hand, since $\sigma_\chi\gg\log q^*$ while $\nu(q,a)\le2^{\omega(q)+1}$, the condition $\Delta_{a,b}>0$ fails for every pair $a\ne b$ once $q$ is large (Proposition~\ref{prop:Delta}).

Rubinstein and Sarnak \cite{RS} showed, under GRH and the linear independence hypothesis LI for the Dirichlet $L$-functions modulo $q$, that the ``prime number race'' between residue classes has limiting logarithmic densities (see also \cite{GM}). The centre-of-mass race has the same densities.

\begin{theorem}\label{thm:race}
Let $q\ge3$ and let $a_1,\dots,a_r$ be distinct residues modulo $q$ coprime to $q$. Assume GRH and LI for the Dirichlet $L$-functions modulo $q$. Then for every permutation $\tau$ of $\{1,\dots,r\}$ the set of $x$ with
\[
D(x;q,a_{\tau(1)})>D(x;q,a_{\tau(2)})>\dots>D(x;q,a_{\tau(r)})
\]
has a logarithmic density, and it is equal to the logarithmic density of the set of $x$ with $\pi(x;q,a_{\tau(1)})>\dots>\pi(x;q,a_{\tau(r)})$.
\end{theorem}

Since $\pi(x;4,3)>\pi(x;4,1)$ on a set of logarithmic density $0.9959\ldots$ \cite{RS}, the centre of mass of the primes $\equiv3\pmod4$ lies further from $\log x$ than that of the primes $\equiv1\pmod4$ on a set of the same logarithmic density.

The centre of mass is the first of a sequence of moments. For $j\ge1$ let
\[
M_j(x)=\frac{1}{\pi(x)}\sum_{p\le x}\Bigl(\log\frac xp\Bigr)^j,\qquad M_j^0(x)=\frac{1}{\Li(x)}\int_2^x\Bigl(\log\frac xt\Bigr)^j\frac{dt}{\log t},
\]
so that $M_1=D$ and $M_1^0=D_0$, and $M_2$ is the moment of inertia of the system about $\log x$. Both $M_j(x)$ and $M_j^0(x)$ tend to $j!$, the $j$-th moment of the exponential law.

\begin{theorem}\label{thm:moments}
Assume the Riemann hypothesis, and let $j\ge1$.
\begin{enumerate}
\item[(a)] $M_j(x)-M_j^0(x)=\Omega_\pm\bigl(x^{-1/2}\log\log\log x\bigr)$.
\item[(b)] The function $\frac{\log x}{j!\sqrt x}\,\pi(x)\bigl(M_j(x)-M_j^0(x)\bigr)$ has a limiting logarithmic distribution with mean $-(2^j-1)$. For $j\ge2$ it therefore differs from that of $\frac{\log x}{\sqrt x}\bigl(\pi(x)-\li(x)\bigr)$; for $j=1$ the two coincide if LI holds.
\item[(c)] If LI also holds, the set of $x$ with $M_j(x)>M_j^0(x)$ has a positive logarithmic density $\delta_j$, and
\[
\log\log\frac1{\delta_j}=\bigl(\sqrt{2\pi}+o(1)\bigr)\,2^{j/2}\qquad(j\to\infty).
\]
\end{enumerate}
\end{theorem}

Thus $\delta_j=\exp\bigl(-\exp((\sqrt{2\pi}+o(1))2^{j/2})\bigr)$. Since the weight of $\rho$ in the $j$-th moment differs from $1/\rho$ by $O(|\rho|^{-2})$, part (c) is in effect the large-deviation estimate of Montgomery \cite{Montgomery} for the variable of Rubinstein and Sarnak, taken at the level $2^j-1$.

The same questions can be asked in the ring $\mathbb F_q[T]$, where the Riemann hypothesis is Weil's theorem. Let $M\in\mathbb F_q[T]$ be monic, and for $A$ prime to $M$ and an integer $X\ge1$ place a unit mass at the point $\deg P$ for every monic irreducible $P\equiv A\pmod M$ with $\deg P\le X$. Let $D_q(X;M,A)$ be the distance of the centre of mass of this system from the point $X$, and let $\pi_q(X;M,A)$ be the number of such $P$. For a non-principal character $\chi$ modulo $M$, write the inverse zeros of $L(u,\chi)$ of modulus $\sqrt q$ as $\sqrt q\,e^{i\alpha_\gamma}$. Following Cha \cite[Definition~3.3]{Cha}, we say that GSH holds if the multiset of the numbers $\alpha_\gamma\in[0,\pi]$ arising from all non-principal characters modulo $M$, counted with multiplicity, is linearly independent over $\mathbb Q$ together with $2\pi$. Thus every such inverse zero is simple and belongs to a single character; this is the multiset form of Cha's set-theoretic condition. Cha \cite{Cha} proved, under GSH, that the race between the $\pi_q(X;M,A)$ has limiting natural densities.

\begin{theorem}\label{thm:FF}
Let $A_1,\dots,A_r$ be distinct residues prime to $M$, and assume GSH for the non-principal characters modulo $M$. Assume further that for all $i\ne j$ some character $\chi$ modulo $M$ with $\chi(A_i)\ne\chi(A_j)$ has an inverse zero of modulus $\sqrt q$. Then for every permutation $\tau$ the set of integers $X$ with
\[
D_q(X;M,A_{\tau(1)})>\dots>D_q(X;M,A_{\tau(r)})
\]
has a natural density, equal to the natural density of the set of $X$ with $\pi_q(X;M,A_{\tau(1)})>\dots>\pi_q(X;M,A_{\tau(r)})$. Moreover, the density of the first set among the even $X$ equals the density of the second set among the odd $X$, and conversely.
\end{theorem}

The hypothesis on the inverse zeros is not automatic: for $\deg M\le2$ some $L(u,\chi)$ has no inverse zero of modulus $\sqrt q$. The factor relating the two explicit formulae has the constant modulus $\sqrt q/(q-1)$, which does not affect the race. The bias terms, on the other hand, are not simply rescaled: for each parity of $X$ the centre-of-mass race carries the bias that the prime race has for the other parity.

The mirror theorems (Theorems~\ref{thm:mirror} and~\ref{thm:mirrorAP}) need the Riemann hypothesis or GRH, but no linear independence hypothesis. The methods are classical: those of Landau and Littlewood for Theorem~\ref{thm:main}, of Rubinstein and Sarnak, in the form given by Akbary, Ng and Shahabi, for Theorems~\ref{thm:skewes} and~\ref{thm:race}, of Montgomery for Theorem~\ref{thm:moments}(c), and of Cha for Theorem~\ref{thm:FF}. What is new is the identity $D=F/\pi$ and the explicit formula for $D-D_0$, from which the criterion of Corollary~\ref{cor:Gn}, the mirror theorems and the interchange of parities are derived.

Section~\ref{sec:mass} contains the mass identities, and Section~\ref{sec:decomp} the decomposition of $D-D_0$ and its explicit formula. Theorem~\ref{thm:main} and Corollary~\ref{cor:Gn} are proved in Section~\ref{sec:main}, Theorems~\ref{thm:mirror}--\ref{thm:race} in Section~\ref{sec:mirror}, Theorem~\ref{thm:moments} in Section~\ref{sec:moments} and Theorem~\ref{thm:FF} in Section~\ref{sec:FF}; a concluding remark concerns other weights.

\section{The mass system}\label{sec:mass}

\begin{lemma}\label{lem:identity}
For $x\ge2$ put $F(x)=\int_2^x\pi(t)\,dt/t$ and $F_0(x)=\int_2^x\Li(t)\,dt/t$. Then
\[
F(x)=\sum_{p\le x}\log\frac xp=\pi(x)\log x-\theta(x),\qquad F_0(x)=\Li(x)\log x-x+2,
\]
and consequently $D(x)=F(x)/\pi(x)$ and $D_0(x)=F_0(x)/\Li(x)$.
\end{lemma}

\begin{proof}
By Stieltjes integration by parts, $\theta(x)=\int_{2^-}^x\log t\,d\pi(t)=\pi(x)\log x-F(x)$. Similarly $\int_2^x\log t\,d\Li(t)=x-2$ gives the formula for $F_0$.
\end{proof}

The function $\log F$ is continuous and piecewise continuously differentiable on $(2,\infty)$, with
\begin{equation}\label{eq:logF}
(\log F)'(x)=\frac{\pi(x)}{xF(x)}=\frac{1}{xD(x)},\qquad(\log F_0)'(x)=\frac{1}{xD_0(x)},
\end{equation}
the first except at the primes. The expansion of $D_0$ follows from Lemma~\ref{lem:identity} and $\Li(x)=\frac{x}{\log x}\sum_{k<K}\frac{k!}{(\log x)^k}+O(x(\log x)^{-K-1})$. Writing $E(x)=\pi(x)-\Li(x)$, Lemma~\ref{lem:identity} also gives
\begin{equation}\label{eq:DD0}
D-D_0=\frac{F-F_0}{\pi}-\frac{F_0\,E}{\pi\,\Li},\qquad F(x)-F_0(x)=\int_2^x\frac{E(t)}{t}\,dt .
\end{equation}

\section{The basic decomposition}\label{sec:decomp}

Let
\[
H(x)=\pi(x)\bigl(D(x)-D_0(x)\bigr),
\]
which has the sign of $D(x)-D_0(x)$.

\begin{lemma}\label{lem:H}
For $x\ge2$,
\begin{align}
H(x)&=-D_0(x)E(x)+\int_2^x\frac{E(t)}{t}\,dt,\label{eq:H1}\\
\int_2^x\frac{E(t)}{t}\,dt&=E(x)\log x-\bigl(\theta(x)-x\bigr)-2,\label{eq:H2}\\
E(x)\log x-\bigl(\theta(x)-x\bigr)&=\log x\int_2^x\frac{\theta(t)-t}{t(\log t)^2}\,dt+O(\log x).\label{eq:H3}
\end{align}
\end{lemma}

\begin{proof}
By Lemma~\ref{lem:identity}, $H=F-D_0\pi=(F_0+\int_2^xE(t)\,dt/t)-D_0(\Li+E)$, and $F_0=D_0\Li$; this is \eqref{eq:H1}. Subtracting the two formulas of Lemma~\ref{lem:identity} gives \eqref{eq:H2}. Finally, partial summation gives $\pi(x)=\frac{\theta(x)}{\log x}+\int_2^x\frac{\theta(t)}{t(\log t)^2}dt$ and $\Li(x)=\frac{x}{\log x}+\int_2^x\frac{dt}{(\log t)^2}+O(1)$, and subtracting and multiplying by $\log x$ gives \eqref{eq:H3}.
\end{proof}

The integral in \eqref{eq:H1} is smoother than $E(x)$; its explicit formula converges absolutely.

\begin{lemma}\label{lem:integral}
Assume the Riemann hypothesis, and write the non-trivial zeros as $\rho=\frac12+i\gamma$. Then
\[
\int_2^x\frac{E(t)}{t}\,dt=-\frac{\sqrt x}{\log x}\Bigl(2+\sum_{\rho}\frac{x^{i\gamma}}{\rho^2}\Bigr)+O\Bigl(\frac{\sqrt x}{(\log x)^2}\Bigr),
\]
where the series converges absolutely. Since $\sigma_0<2$, it follows that $\int_2^xE(t)\,dt/t<0$ for all sufficiently large $x$.
\end{lemma}

\begin{proof}
By \eqref{eq:H2} and \eqref{eq:H3} it suffices to show that
\[
\int_2^x\frac{\theta(t)-t}{t(\log t)^2}\,dt=-\frac{\sqrt x}{(\log x)^2}\Bigl(2+\sum_\rho\frac{x^{i\gamma}}{\rho^2}\Bigr)+O\Bigl(\frac{\sqrt x}{(\log x)^3}\Bigr).
\]
Write $\theta(t)-t=(\psi(t)-t)-(\psi(t)-\theta(t))$. Under the Riemann hypothesis $\psi(t)-\theta(t)=\theta(t^{1/2})+O(t^{1/3}\log t)=\sqrt t+O(t^{1/3}\log t)$. The contribution of $-\sqrt t$ is $-\int_2^x\frac{dt}{\sqrt t(\log t)^2}=-\frac{2\sqrt x}{(\log x)^2}+O\bigl(\frac{\sqrt x}{(\log x)^3}\bigr)$, and that of the term $O(t^{1/3}\log t)$ is $O(x^{1/3})$. For the contribution of $\psi(t)-t$ we use the explicit formula
\[
\psi_1(t):=\int_0^t\bigl(\psi(u)-u\bigr)\,du=-\sum_\rho\frac{t^{\rho+1}}{\rho(\rho+1)}+O(t),
\]
valid for $t\ge2$ (see \cite[Ch.~IV]{Ingham} or \cite[\S12.1]{MV}). Under the Riemann hypothesis $|t^{\rho+1}|=t^{3/2}$, so the series converges absolutely and uniformly on $[2,x]$, since $\sum_\rho|\rho(\rho+1)|^{-1}<\infty$. Integrating by parts with $g(t)=1/(t(\log t)^2)$ gives
\[
\int_2^x(\psi(t)-t)g(t)\,dt=\psi_1(x)g(x)-\int_2^x\psi_1(t)g'(t)\,dt+O(1).
\]
Inserting the explicit formula and integrating term by term, which is legitimate by uniform convergence, the zero $\rho$ contributes (the boundary term at $t=2$ being $O(|\rho|^{-2})$)
\[
\begin{aligned}
-\frac{1}{\rho(\rho+1)}\Bigl(x^{\rho+1}g(x)-\int_2^xt^{\rho+1}g'(t)\,dt\Bigr)&=-\frac1\rho\int_2^x\frac{t^{\rho-1}}{(\log t)^2}\,dt+O\Bigl(\frac{1}{|\rho|^2}\Bigr)\\
&=-\frac{x^{\rho}}{\rho^2(\log x)^2}+O\Bigl(\frac{\sqrt x}{|\rho|^2(\log x)^3}\Bigr),
\end{aligned}
\]
the last step by one more integration by parts. Since $|t^\rho|=\sqrt t$, the implied constants are uniform in $\rho$. Summing over $\rho$ with $\sum_\rho|\rho|^{-2}<\infty$, and noting that the term $O(t)$ in $\psi_1$ contributes $O\bigl(xg(x)+\int_2^xt|g'(t)|\,dt\bigr)=O(1)$, we obtain the assertion. Under the Riemann hypothesis $|\rho|^2=\rho(1-\rho)$, so the sum $\sum_\rho|\rho|^{-2}$ is the constant $\sigma_0$ of the introduction.
\end{proof}

\begin{remark}
Lemma~\ref{lem:integral} is the analogue for the primes of the classical explicit formula
\[
\sum_{n\le x}\Lambda(n)\log\frac xn=\int_1^x\psi(t)\,\frac{dt}{t}=x-\sum_\rho\frac{x^\rho}{\rho^2}+O(\log x)
\]
for the Riesz mean; the constant $-2$, which has no counterpart in the Riesz mean of $\psi$, comes from the squares of primes. Under the Riemann hypothesis it shows that the moment $F(x)$ eventually stays below its model $F_0(x)$; the converse is Corollary~\ref{cor:Riesz} below.
\end{remark}

\section{Proof of Theorem~\ref{thm:main}}\label{sec:main}

For $\operatorname{Re}s>1$ we have
\[
\int_2^\infty\bigl(F(x)-F_0(x)\bigr)x^{-s-1}\,dx=\frac{1}{s^2}\Bigl(P(s)-\int_2^\infty\frac{x^{-s}}{\log x}\,dx\Bigr)=:\frac{\mathcal M(s)}{s^2},
\]
where $P(s)=\sum_pp^{-s}$; this follows from $F'=\pi(x)/x$, $F_0'=\Li(x)/x$ and $\int_2^\infty\pi(x)x^{-s-1}dx=P(s)/s$. Now $P(s)=\sum_{k\ge1}\frac{\mu(k)}{k}\log\zeta(ks)$, so $P(s)-\log\zeta(s)$ is analytic for $\operatorname{Re}s>\frac12$. The function $\log\zeta(s)+\log(s-1)$ is analytic in a neighbourhood of the real half-line $s>\frac12$, because $\zeta$ has no zeros on the real segment $(0,1)$; and $\int_2^\infty x^{-s}dx/\log x+\log(s-1)$ is entire. Hence $\mathcal M(s)$ is analytic in a neighbourhood of the real half-line $s>\frac12$, while it has a logarithmic singularity at every zero $\rho$ of $\zeta$ with $\operatorname{Re}\rho>\frac12$.

\begin{proof}[Proof of Theorem~\ref{thm:main}(a)]
By a classical theorem (see \cite[Ch.~IV]{Ingham}), $E(x)\ll_\varepsilon x^{\Theta+\varepsilon}$, so $F-F_0\ll x^{\Theta+\varepsilon}$ and \eqref{eq:DD0} gives $D-D_0\ll x^{\Theta-1+2\varepsilon}$.

Now let $\Theta>\frac12$ and $\frac12\le\beta<\Theta$; it suffices to treat such $\beta$. Suppose that, for some $K>0$, $D(x)-D_0(x)\le Kx^{\beta-1}$ for all $x\ge x_0$. Since $D$ and $D_0$ tend to $1$, \eqref{eq:logF} gives, except at the primes,
\[
\bigl(\log F-\log F_0\bigr)'(x)=\frac{D_0(x)-D(x)}{x\,D(x)\,D_0(x)}\ge-2Kx^{\beta-2}\qquad(x\ge x_1).
\]
Integrating over $[x,y]$, letting $y\to\infty$ and using $\log(F/F_0)\to0$, we obtain $\log(F/F_0)(x)\le\frac{2K}{1-\beta}x^{\beta-1}$. Since $F_0(x)\ll x/\log x$, this gives $F(x)-F_0(x)\le Cx^\beta$ for $x\ge x_1$. Hence $g(x)=Cx^\beta-(F(x)-F_0(x))$ is non-negative for $x\ge x_1$, and for $\operatorname{Re}s>1$
\[
\int_2^\infty g(x)x^{-s-1}\,dx=\frac{C\,2^{\beta-s}}{s-\beta}-\frac{\mathcal M(s)}{s^2}.
\]
The right-hand side is analytic in a neighbourhood of the real half-line $s>\beta$. As the integral over $[2,x_1]$ is entire, Landau's theorem on integrals with non-negative integrands (see \cite[Ch.~V]{Ingham}) shows that the integral converges, and represents an analytic function, for $\operatorname{Re}s>\beta$. This is impossible, since $\mathcal M(s)$ has a singularity at a zero $\rho$ with $\operatorname{Re}\rho>\beta$. Hence $\limsup(D-D_0)x^{1-\beta}=+\infty$. If instead $D-D_0\ge-Kx^{\beta-1}$ for $x\ge x_0$, the same argument applied to $Cx^\beta+(F-F_0)$ gives a contradiction, so the $\liminf$ is $-\infty$.
\end{proof}

\begin{proof}[Proof of Theorem~\ref{thm:main}(b)]
If the Riemann hypothesis fails, (b) follows from (a) with $\frac12<\beta<\Theta$; we may therefore assume the Riemann hypothesis. By \eqref{eq:H1} and Lemma~\ref{lem:integral}, $H(x)=-D_0(x)E(x)+O(\sqrt x/\log x)$, where $D_0(x)\to1$. Since $\Li$ and $\li$ differ by a constant, Littlewood's theorem (see \cite[Ch.~V]{Ingham}) gives
\[
E(x)=\Omega_\pm\Bigl(\frac{\sqrt x\log\log\log x}{\log x}\Bigr).
\] At points where $E(x)>c\sqrt x\log\log\log x/\log x$ we get $H(x)<-\frac c2\sqrt x\log\log\log x/\log x$ for large $x$, and similarly with the signs reversed. Dividing by $\pi(x)\sim x/\log x$ gives the assertion.
\end{proof}

The remaining assertions follow at once. If the Riemann hypothesis holds, (b) gives the sign changes and shows that $D-D_0\ll x^\alpha$ fails for $\alpha\le-\frac12$; otherwise $\Theta>\frac12$ and (a) gives both. Finally, if $D-D_0\ll x^{-1/2+\varepsilon}$ for every $\varepsilon>0$, then $\Theta-1\le-\frac12$, which is the Riemann hypothesis.

\begin{proof}[Proof of Corollary~\ref{cor:Gn}]
(a) If the Riemann hypothesis holds, the assertion is Theorem~\ref{thm:main} at $x=p_n$. Conversely, put $d_n=D(p_n)=\log p_n-\theta(p_n)/n$ and assume that $d_n=D_0(p_n)+O(p_n^{-1/2+\varepsilon})$. Since $\theta(p_{n+1})=\theta(p_n)+\log p_{n+1}$, the number $g_n=\log(p_{n+1}/p_n)$ satisfies
\[
g_n=\frac{n+1}{n}\Bigl(d_{n+1}-d_n+\frac{d_n}{n+1}\Bigr).
\]
Since $D_0'(x)=x^{-1}\bigl(1-xD_0(x)/(\Li(x)\log x)\bigr)$, the expansions of $D_0$ and $\Li$ give $|D_0'(x)|\le C/(x\log^2x)$; hence $|D_0(p_{n+1})-D_0(p_n)|\le Cg_n/\log^2p_n$, so that $d_{n+1}-d_n=O(g_n/\log^2p_n)+O(p_n^{-1/2+\varepsilon})$; and $d_n/n\ll\log p_n/p_n$. Hence $g_n\ll p_n^{-1/2+\varepsilon}$. For $p_n\le x<p_{n+1}$ the functions $\theta$ and $\pi$ are constant on $[p_n,x]$, so $D(x)=d_n+\log(x/p_n)$, while $D_0(x)-D_0(p_n)\ll\log(x/p_n)$; and $p_n^{-1/2+\varepsilon}\ll x^{-1/2+\varepsilon}$ since $p_{n+1}<2p_n$. Therefore $D(x)=D_0(x)+O(x^{-1/2+\varepsilon})$ for all $x\ge2$, and Theorem~\ref{thm:main} gives the Riemann hypothesis.

(b) On $[p_n,p_{n+1})$ we have $D'(x)=1/x$ and $|D_0'(x)|\le C/(x\log^2x)$, so $D-D_0$ is increasing there for large $n$. At $p_{n+1}$ the function $D$ jumps by $\theta(p_n)/n-\theta(p_{n+1})/(n+1)=-(d_n+g_n)/(n+1)$, which is negative and $O(\log p_n/p_n)$, while $D_0$ is continuous. Hence, for $p_n\le x<p_{n+1}$,
\[
d_n-D_0(p_n)\le D(x)-D_0(x)\le d_{n+1}-D_0(p_{n+1})+O\Bigl(\frac{\log p_n}{p_n}\Bigr).
\]
By Theorem~\ref{thm:main}(b) there are arbitrarily large $x$ with $D(x)-D_0(x)>cx^{-1/2}\log\log\log x$, and for these the right-hand inequality gives $d_{n+1}-D_0(p_{n+1})>\frac c2p_{n+1}^{-1/2}\log\log\log p_{n+1}$. Similarly, there are arbitrarily large $x$ with $D(x)-D_0(x)<-cx^{-1/2}\log\log\log x$, and for these the left-hand inequality and $p_{n+1}<2p_n$ give $d_n-D_0(p_n)<-\frac c2p_n^{-1/2}\log\log\log p_n$.
\end{proof}

The argument of Theorem~\ref{thm:main}(a) also gives a criterion for the Riemann hypothesis in terms of the moment $F$ alone.

\begin{corollary}\label{cor:Riesz}
The Riemann hypothesis holds if and only if $\int_2^x(\pi(t)-\Li(t))\,dt/t<0$ for all sufficiently large $x$ (compare \cite[Theorem~1.4]{Johnston}).
\end{corollary}

\begin{proof}
If the Riemann hypothesis holds, the integral is negative for large $x$ by Lemma~\ref{lem:integral}. Conversely, suppose that $F(x)-F_0(x)=\int_2^xE(t)\,dt/t\le0$ for $x\ge x_1$. Then $g=F_0-F$ is non-negative for $x\ge x_1$, and $\int_2^\infty g(x)x^{-s-1}\,dx=-\mathcal M(s)/s^2$ for $\operatorname{Re}s>1$. The right-hand side is analytic in a neighbourhood of the real half-line $s>\frac12$, and the integral over $[2,x_1]$ is entire, so Landau's theorem shows that the integral converges for $\operatorname{Re}s>\frac12$. Hence $\mathcal M(s)$ is analytic there, and $\zeta$ has no zeros with real part greater than $\frac12$.
\end{proof}

The converse in Corollary~\ref{cor:Riesz} also follows from the case $c=1$ of a theorem of Johnston \cite[Theorem~1.4]{Johnston}, which shows that if $\Theta>\frac12$ and $c<1+\Theta$, then $\int_2^x(\pi(t)-\li(t))t^{-c}\,dt$ takes positive values for arbitrarily large $x$; replacing $\li$ by $\Li$ adds $\li(2)\log(x/2)>0$. For the unweighted integral, Pintz \cite{Pintz} asserted that the Riemann hypothesis is equivalent to $\int_2^x(\pi(t)-\li(t))\,dt<0$ for large $x$, and Johnston \cite{Johnston} gave a proof, with the inequality for all $x>2$.

\section{Limiting distributions and the mirror theorems}\label{sec:mirror}

We first prove Theorem~\ref{thm:mirror}, which needs only \eqref{eq:H1} and Lemma~\ref{lem:integral}.

\begin{proof}[Proof of Theorem~\ref{thm:mirror}]
Write $\widehat E=(\log x/\sqrt x)(\pi(x)-\Li(x))$; replacing $\Li$ by $\li$ changes $\widehat E$ by $O(\log x/\sqrt x)$. By \eqref{eq:H1} and Lemma~\ref{lem:integral},
\[
\frac{\log x}{\sqrt x}H(x)=-D_0(x)\widehat E-2-\Sigma(x)+O\Bigl(\frac1{\log x}\Bigr),
\]
where $\Sigma(x)=\sum_\rho x^{i\gamma}/\rho^2$ satisfies $|\Sigma(x)|\le\sigma_0$.
If $\widehat E<-(2+\sigma_0+\varepsilon)$, then $-D_0\widehat E>2+\sigma_0+\varepsilon$, since $D_0(x)>1$ for large $x$, so $H(x)>0$ for large $x$. If $H(x)>0$, then $D_0\widehat E<-2+\sigma_0+O(1/\log x)$, and since $D_0(x)=1+O(1/\log x)$ we get $\widehat E<-(2-\sigma_0-\varepsilon)$ for large $x$.
\end{proof}

We use the method of Rubinstein and Sarnak \cite{RS}, in the form of $B^2$-almost periodic functions described by Akbary, Ng and Shahabi \cite[\S2]{ANS}. Let $\Psi$ be a real function of $y\ge y_0$, and suppose that, for every $T\ge2$,
\begin{equation}\label{eq:B2}
\Psi(y)=c+\operatorname{Re}\sum_{0<\gamma\le T}r_\gamma e^{i\gamma y}+\varepsilon_T(y),\qquad\lim_{T\to\infty}\limsup_{Y\to\infty}\frac1Y\int_{y_0}^Y|\varepsilon_T(y)|^2\,dy=0,
\end{equation}
where $\gamma$ runs over an increasing sequence of positive numbers and $\sum_\gamma|r_\gamma|^2<\infty$. Then $\Psi$ is $B^2$-almost periodic, and therefore has a limiting distribution $\mu$, in the sense that
\[
\frac1Y\int_{y_0}^Yf(\Psi(y))\,dy\to\int f\,d\mu
\]
for bounded continuous $f$ (Besicovitch; see \cite[Theorem~2.9]{ANS}). If, in addition, the $\gamma$ are linearly independent over $\mathbb Q$, then $\mu$ is the distribution of the random variable $c+\operatorname{Re}\sum_\gamma r_\gamma Z_\gamma$, where the $Z_\gamma$ are independent and uniformly distributed on the unit circle; see \cite[\S3]{RS} and \cite[Theorem~1.9]{ANS}; the functions considered below have $|r_\gamma|\ll1/\gamma$ and the mean-square bounds of Rubinstein and Sarnak, as Lemmas~\ref{lem:normal} and~\ref{lem:normalAP} show. In that case $\mu$ is absolutely continuous as soon as one $r_\gamma$ is non-zero, since $\operatorname{Re}(r_\gamma Z_\gamma)$ then has the arcsine law; so the logarithmic density of $\{x:\Psi(\log x)>0\}$ exists and equals $\mu((0,\infty))$. The same holds for vector-valued functions, with the joint distribution of the corresponding random vector \cite{RS,ANS}.

The law of $c+\operatorname{Re}\sum_\gamma r_\gamma Z_\gamma$ depends only on $c$ and on the moduli $|r_\gamma|$, since $Z_\gamma$ and $u_\gamma Z_\gamma$ have the same law for every $u_\gamma$ of modulus one. In the vector case it suffices that the factor $u_\gamma$ depends only on $\gamma$ and not on the coordinate.

\begin{lemma}\label{lem:normal}
Assume the Riemann hypothesis, and put $y=\log x$. For every $T\ge2$,
\[
\frac{\log x}{\sqrt x}H(x)=-1+2\operatorname{Re}\sum_{0<\gamma\le T}x^{i\gamma}\Bigl(\frac1\rho-\frac1{\rho^2}\Bigr)+\varepsilon_T(y),
\]
where
\[
\limsup_{Y\to\infty}\frac1Y\int_1^Y|\varepsilon_T(y)|^2\,dy\ll\frac{(\log T)^2}{T}.
\]
\end{lemma}

\begin{proof}
Write $\widehat E(y)=\frac{\log x}{\sqrt x}E(x)$. By \eqref{eq:H1},
\[
\frac{\log x}{\sqrt x}H(x)=-\widehat E(y)-\bigl(D_0(x)-1\bigr)\widehat E(y)+\frac{\log x}{\sqrt x}\int_2^x\frac{E(t)}{t}\,dt ,
\]
and we treat the three terms in turn. First, by the explicit formula for $\pi(x)$ and the mean-square estimate of Rubinstein and Sarnak \cite[\S2]{RS} (see also \cite{ANS}),
\[
\begin{gathered}
\widehat E(y)=-1-2\operatorname{Re}\sum_{0<\gamma\le T}\frac{x^{i\gamma}}{\rho}+\varepsilon^{(1)}_T(y),\\
\limsup_{Y\to\infty}\frac1Y\int_1^Y|\varepsilon^{(1)}_T|^2\,dy\ll\frac{(\log T)^2}{T};
\end{gathered}
\]
the bias $-1$ comes from the squares of primes, and replacing $\li$ by $\Li$ changes $\widehat E$ by $O(y\,e^{-y/2})$. Secondly, Cram\'er's mean-value theorem \cite{Cramer} states that $\frac1Y\int_1^Y|(\psi(e^y)-e^y)e^{-y/2}|^2dy\ll1$ under the Riemann hypothesis; by \eqref{eq:H3}, together with the bound $\int_2^x(\theta(t)-t)\,dt/(t\log^2t)\ll\sqrt x/(\log x)^2$ obtained in the proof of Lemma~\ref{lem:integral}, it follows that $\frac1Y\int_1^Y|\widehat E(y)|^2dy\ll1$, and $D_0(x)-1\ll1/y$, so the mean square of the second term tends to $0$. Thirdly, by Lemma~\ref{lem:integral} the last term equals
\[
-2-2\operatorname{Re}\sum_{0<\gamma\le T}\frac{x^{i\gamma}}{\rho^2}+O\Bigl(\sum_{|\gamma|>T}\frac1{|\rho|^2}+\frac1y\Bigr)=-2-2\operatorname{Re}\sum_{0<\gamma\le T}\frac{x^{i\gamma}}{\rho^2}+O\Bigl(\frac{\log T}{T}+\frac1y\Bigr),
\]
uniformly in $x$. Adding the three contributions gives the lemma.
\end{proof}

\begin{proof}[Proof of Theorem~\ref{thm:skewes}]
By Lemma~\ref{lem:normal}, the function $\Psi(y)=\frac{\log x}{\sqrt x}H(x)$ satisfies \eqref{eq:B2} with $c=-1$ and $r_\gamma=2(1/\rho-1/\rho^2)$. By \cite[\S2]{RS}, the function $\frac{\log x}{\sqrt x}(\pi(x)-\li(x))$ satisfies \eqref{eq:B2} with $c=-1$ and $r_\gamma=-2/\rho$. For $\rho=\frac12+i\gamma$ we have $\rho-1=-\bar\rho$, so
\[
\frac1\rho-\frac1{\rho^2}=\frac{\rho-1}{\rho^2}=u_\gamma\cdot\frac{1}{\rho},\qquad u_\gamma=-\frac{\bar\rho}{\rho},\quad|u_\gamma|=1 .
\]
Under LI, therefore, both functions have the limiting distribution of $-1+\operatorname{Re}\sum_{\gamma>0}\frac{2}{|\rho|}Z_\gamma$, which is absolutely continuous. Since $\{x:D(x)>D_0(x)\}=\{x:H(x)>0\}$, the logarithmic densities of this set and of $\{x:\pi(x)>\li(x)\}$ exist and are equal.
\end{proof}

We now turn to arithmetic progressions. Fix $q\ge3$ and write $\varphi=\varphi(q)$. For $(a,q)=1$ put
\[
E(x;q,a)=\pi(x;q,a)-\frac{\Li(x)}{\varphi},\quad H(x;q,a)=\pi(x;q,a)\bigl(D(x;q,a)-D_0(x)\bigr).
\]
GRH and LI are as in the introduction. In Lemma~\ref{lem:normalAP} and Theorem~\ref{thm:race} we also assume that $L(\frac12,\chi)\ne0$ for every $\chi\ne\chi_0$ modulo $q$, which is implied by LI; recall that $\zeta(\frac12)\ne0$. For a character $\chi$ modulo $q$ we write the non-trivial zeros of $L(s,\chi)$ (of the primitive character inducing $\chi$) as $\rho_\chi=\frac12+i\gamma_\chi$.

\begin{lemma}\label{lem:integralAP}
Assume GRH. Then
\[
\frac{\varphi\log x}{\sqrt x}\int_2^x\frac{E(t;q,a)}{t}\,dt=-2\nu(q,a)-\sum_\chi\bar\chi(a)\sum_{\rho_\chi}\frac{x^{i\gamma_\chi}}{\rho_\chi^2}+O_q\Bigl(\frac{1}{\log x}\Bigr),
\]
where $\chi$ runs over all characters modulo $q$ and the double series converges absolutely.
\end{lemma}

\begin{proof}
The proofs of Lemmas~\ref{lem:identity} and~\ref{lem:H} apply without change to the primes $p\equiv a\pmod q$. With $F(x;q,a)=\int_2^x\pi(t;q,a)\,dt/t=\pi(x;q,a)\log x-\theta(x;q,a)$ and $F_0/\varphi$ in place of $F$ and $F_0$, they give
\[
\int_2^x\frac{E(t;q,a)}{t}\,dt=\log x\int_2^x\frac{\theta(t;q,a)-t/\varphi}{t(\log t)^2}\,dt+O(\log x).
\]
By orthogonality, $\theta(t;q,a)-t/\varphi=\frac1\varphi\sum_\chi\bar\chi(a)\bigl(\psi(t,\chi)-e_\chi t\bigr)-\bigl(\psi(t;q,a)-\theta(t;q,a)\bigr)$, where $e_\chi=1$ if $\chi$ is principal and $0$ otherwise. Under GRH the prime number theorem for the progressions modulo $q$ with error $O(t^{1/2+\varepsilon})$ gives
\[
\psi(t;q,a)-\theta(t;q,a)=\sum_{\substack{p^2\le t\\ p^2\equiv a}}\log p+O(t^{1/3}\log t)=\frac{\nu(q,a)}{\varphi}\sqrt t+O_q(t^{1/3}\log t),
\]
whose contribution is $-\frac{2\nu(q,a)}{\varphi}\frac{\sqrt x}{(\log x)^2}+O_q\bigl(\frac{\sqrt x}{(\log x)^3}\bigr)$, exactly as in Lemma~\ref{lem:integral}. For each $\chi$ we use the explicit formula
\[
\int_0^t\bigl(\psi(u,\chi)-e_\chi u\bigr)\,du=-\sum_{\rho_\chi}\frac{t^{\rho_\chi+1}}{\rho_\chi(\rho_\chi+1)}+O_q(t\log t)\qquad(t\ge2),
\]
which follows from the explicit formula for $\psi(t,\chi)$ \cite[Ch.~19]{Davenport} by integration, or directly by Mellin inversion \cite[\S12.1]{MV}; replacing $\chi$ by the primitive character inducing it changes $\psi(t,\chi)$ by $O_q(\log t)$. Under GRH the series converges absolutely and uniformly on $[2,x]$. The argument of Lemma~\ref{lem:integral}, with $g(t)=1/(t(\log t)^2)$, now gives, with implied constants uniform in $\rho_\chi$, the contribution $-\frac1\varphi\bar\chi(a)\sum_{\rho_\chi}\frac{x^{\rho_\chi}}{\rho_\chi^2(\log x)^2}+O_q\bigl(\frac{\sqrt x}{(\log x)^3}\bigr)$ from each $\chi$; the term $O_q(t\log t)$ contributes $O_q(\log x)$. Multiplying by $\log x$ and normalizing gives the lemma.
\end{proof}

\begin{lemma}\label{lem:HtoD}
Assume GRH, and let $a,b$ be residues prime to $q$. Then
\[
D(x;q,a)-D(x;q,b)=\frac{\varphi}{\Li(x)}\bigl(H(x;q,a)-H(x;q,b)\bigr)+O\Bigl(\frac{(\log x)^6}{x}\Bigr).
\]
\end{lemma}

\begin{proof}
Under GRH, $E(x;q,a)\ll\sqrt x(\log x)^2$ and $H(x;q,a)\ll\sqrt x(\log x)^2$. Writing $\pi(x;q,a)=\frac{\Li(x)}{\varphi}(1+\eta_a)$ with $\eta_a\ll x^{-1/2}(\log x)^3$, we have $D(x;q,a)-D_0(x)=\frac{\varphi}{\Li(x)}H(x;q,a)(1+\eta_a)^{-1}$, and the lemma follows.
\end{proof}

\begin{proof}[Proof of Theorem~\ref{thm:mirrorAP}]
Put $h_a=\frac{\varphi\log x}{\sqrt x}H(x;q,a)$ and $\widehat E_a=\frac{\varphi\log x}{\sqrt x}E(x;q,a)$. By the analogue of \eqref{eq:H1} and Lemma~\ref{lem:integralAP},
\[
h_a=-D_0\widehat E_a-2\nu(q,a)-\sum_\chi\bar\chi(a)\Sigma_\chi(x)+O_q\Bigl(\frac1{\log x}\Bigr),
\]
where $\Sigma_\chi(x)=\sum_{\rho_\chi}x^{i\gamma_\chi}/\rho_\chi^2$ satisfies $|\Sigma_\chi(x)|\le\sigma_\chi$.
The principal character contributes equally to $h_a$ and $h_b$, so
\[
h_b-h_a\ge-D_0\bigl(\widehat E_b-\widehat E_a\bigr)+\Delta_{a,b}-O_q\Bigl(\frac1{\log x}\Bigr).
\]
By Lemma~\ref{lem:HtoD}, $D(x;q,a)-D(x;q,b)$ has the sign of $h_a-h_b+O((\log x)^6/\sqrt x)$. Hence $D(x;q,a)\ge D(x;q,b)$ implies $h_b-h_a\le o(1)$, and so $D_0(\widehat E_b-\widehat E_a)\ge\Delta_{a,b}-o(1)$. Since $\widehat E_b-\widehat E_a=\frac{\varphi\log x}{\sqrt x}(\pi(x;q,b)-\pi(x;q,a))$ and $D_0\to1$, the assertion follows.
\end{proof}

\begin{lemma}\label{lem:sigma}
Let $\chi$ be a primitive character modulo $q^*>1$ with $\chi(-1)=(-1)^{\mathfrak a}$, $\mathfrak a\in\{0,1\}$, and assume the Riemann hypothesis for $L(s,\chi)$. Then
\[
\sum_{\rho}\frac{1}{|\rho|^2}=2\operatorname{Re}\frac{L'}{L}(1,\chi)+\log\frac{q^*}{\pi}+\frac{\Gamma'}{\Gamma}\Bigl(\frac{1+\mathfrak a}2\Bigr),
\]
the sum running over the non-trivial zeros of $L(s,\chi)$.
\end{lemma}

\begin{proof}
For $\rho=\frac12+i\gamma$ we have $\operatorname{Re}\frac1\rho=\operatorname{Re}\frac{1}{1-\rho}=\frac{1}{2|\rho|^2}$. By the Hadamard product \cite[Ch.~12]{Davenport},
\[
\frac{L'}{L}(s,\chi)=-\frac12\log\frac{q^*}{\pi}-\frac12\frac{\Gamma'}{\Gamma}\Bigl(\frac{s+\mathfrak a}2\Bigr)+b(\chi)+\sum_\rho\Bigl(\frac1{s-\rho}+\frac1\rho\Bigr),
\]
with $\operatorname{Re}b(\chi)=-\sum_\rho\operatorname{Re}\frac1\rho$. Taking real parts at $s=1$ gives the lemma.
\end{proof}

For an odd primitive real character $\chi$ modulo $q^*$, Lemma~\ref{lem:sigma} and the functional equation of $L(s,\chi)$ give
\[
\sigma_\chi=\gamma_0+2\log2-\log\frac{q^*}{\pi}-2\,\frac{L'}{L}(0,\chi).
\]
Since $L(s,\chi)=q^{*-s}\sum_{a=1}^{q^*}\chi(a)\zeta(s,a/q^*)$, Lerch's formula $\zeta'(0,\alpha)=\log\Gamma(\alpha)-\frac12\log2\pi$ for the Hurwitz zeta function \cite{Berndt} gives
\[
\begin{aligned}
\frac{L'}{L}(0,\chi_{-4})&=4\log\Gamma\bigl(\tfrac14\bigr)-2\log\pi-3\log2,\\
\frac{L'}{L}(0,\chi_{-3})&=6\log\Gamma\bigl(\tfrac13\bigr)+\tfrac12\log3-3\log2-3\log\pi,
\end{aligned}
\]
and hence the closed forms of $\sigma_{\chi_{-4}}$ and $\sigma_{\chi_{-3}}$ stated in the introduction. The only non-principal character modulo $6$ is induced by $\chi_{-3}$, and $\nu(6,1)=2$, $\nu(6,5)=0$; this gives the value of $\Delta_{1,5}$ for $q=6$.

\begin{proposition}\label{prop:Delta}
There is $q_0$ such that $\Delta_{a,b}<0$ for every $q\ge q_0$ and all distinct residues $a,b$ prime to $q$.
\end{proposition}

\begin{proof}
Let $\chi$ be primitive modulo $q^*$. The number $N(T,\chi)$ of zeros of $L(s,\chi)$ with $0<\beta<1$ and $|\gamma|\le T$ satisfies $N(T,\chi)=\frac T\pi\log\frac{q^*T}{2\pi e}+O(\log q^*T)$ for $T\ge2$ \cite[Ch.~16]{Davenport}. Fix $T_0\ge2\pi e$ so large that the error term is at most $\frac{T_0}{2\pi}\log q^*$ for all $q^*\ge3$; then $N(T_0,\chi)\ge\frac{T_0}{2\pi}\log q^*$ for all $q^*\ge3$, and since $|\rho_\chi|^2\le1+T_0^2$ for these zeros, $\sigma_\chi\ge c\log q^*$ with $c=T_0/(2\pi(1+T_0^2))$. Moreover $\sigma_\chi\ge0$ for every $\chi$. Since $|\chi(a)-\chi(b)|\le2$, orthogonality gives
\[
\sum_\chi|\chi(a)-\chi(b)|\ge\frac12\sum_\chi|\chi(a)-\chi(b)|^2=\varphi(q)\qquad(a\ne b).
\]
At most $\sum_{d\mid q,\,d\le\sqrt q}\varphi(d)\le\tau(q)\sqrt q$ characters modulo $q$ have conductor $q^*\le\sqrt q$, and they contribute at most $2\tau(q)\sqrt q$ to the left-hand side. For the other characters, if $q\ge9$, we have $\sigma_\chi\ge\frac c2\log q$. Hence
\[
\sum_{\chi\ne\chi_0}|\chi(a)-\chi(b)|\,\sigma_\chi\ge\bigl(\varphi(q)-2\tau(q)\sqrt q\bigr)\,\frac c2\log q,
\]
which exceeds $2\nu(q,a)\le2^{\omega(q)+2}=q^{o(1)}$ for large $q$, since $\varphi(q)\gg q/\log\log q$.
\end{proof}

Proposition~\ref{prop:Delta} is unconditional. It concerns only the sufficient condition $\Delta_{a,b}>0$ of Theorem~\ref{thm:mirrorAP}; it is consistent with the fact that the biases in the prime number race dissolve as $q\to\infty$ \cite{RS,FM}.

\begin{lemma}\label{lem:normalAP}
Assume GRH and $L(\frac12,\chi)\ne0$ for all $\chi\ne\chi_0$ modulo $q$, put $y=\log x$, and let $(a,q)=1$. For every $T\ge2$,
\[
\frac{\varphi\log x}{\sqrt x}H(x;q,a)=-\nu(q,a)+2\operatorname{Re}\sum_\chi\bar\chi(a)\sum_{0<\gamma_\chi\le T}x^{i\gamma_\chi}\Bigl(\frac{1}{\rho_\chi}-\frac1{\rho_\chi^2}\Bigr)+\varepsilon_{a,T}(y),
\]
where $\chi$ runs over all characters modulo $q$ and, for each $\chi$, $\gamma_\chi$ runs over the positive ordinates of the zeros of $L(s,\chi)$; the zeros with negative ordinates are accounted for by the complex conjugate terms, since the zeros of $L(s,\bar\chi)$ are the conjugates of those of $L(s,\chi)$. The errors satisfy
\[
\limsup_{Y\to\infty}\frac1Y\int_1^Y|\varepsilon_{a,T}(y)|^2\,dy\ll_q\frac{(\log T)^2}{T}.
\]
\end{lemma}

\begin{proof}
Put $\widehat E_a(y)=\frac{\varphi\log x}{\sqrt x}E(x;q,a)$. As in Lemma~\ref{lem:normal},
\[
\frac{\varphi\log x}{\sqrt x}H(x;q,a)=-\widehat E_a(y)-\bigl(D_0(x)-1\bigr)\widehat E_a(y)+\frac{\varphi\log x}{\sqrt x}\int_2^x\frac{E(t;q,a)}{t}\,dt.
\]
This is proved as \eqref{eq:H1}, since $D_0\Li/\varphi=F_0/\varphi$. For the first term, $\varphi E(x;q,a)=\bigl(\varphi\pi(x;q,a)-\pi(x)\bigr)+\bigl(\pi(x)-\Li(x)\bigr)$. The two brackets have the expansions of Rubinstein and Sarnak \cite[\S2]{RS}, with biases $-c(q,a)=1-\nu(q,a)$ and $-1$ and with mean-square errors $\ll_q(\log T)^2/T$. Adding them,
\[
\widehat E_a(y)=-\nu(q,a)-2\operatorname{Re}\sum_\chi\bar\chi(a)\sum_{0<\gamma_\chi\le T}\frac{x^{i\gamma_\chi}}{\rho_\chi}+\varepsilon^{(1)}_{a,T}(y).
\]
For the second term, the analogue of Cram\'er's theorem for $\psi(x,\chi)$ under GRH \cite[\S2]{RS} gives $\frac1Y\int_1^Y|\widehat E_a|^2dy\ll_q1$; since $D_0(x)-1\ll1/y$, the mean square of this term tends to $0$. The third term is given by Lemma~\ref{lem:integralAP}, and truncating its absolutely convergent series at $T$ costs $O_q((\log T)/T)$ uniformly. Adding the three contributions proves the lemma.
\end{proof}

\begin{proof}[Proof of Theorem~\ref{thm:race}]
Put $\Psi_j(y)=\frac{\varphi\log x}{\sqrt x}H(x;q,a_j)$ and $\pi_j(y)=\frac{\varphi\log x}{\sqrt x}\bigl(\pi(x;q,a_j)-\pi(x)/\varphi\bigr)$ for $1\le j\le r$. The zeros of $\zeta(s)$ contribute the same terms to every $\Psi_j$, and we remove them by putting $\Psi^*_j=\Psi_j-\frac{\log x}{\sqrt x}H(x)$; this changes no difference $\Psi_i-\Psi_j$, and below we work with the $\Psi^*_j$. The functions $\pi_j$ contain no such terms. Then, by Lemma~\ref{lem:normalAP} and \cite[\S2]{RS}, the vector functions $(\Psi^*_j)_j$ and $(\pi_j)_j$ satisfy the vector form of \eqref{eq:B2}, with the frequencies $\gamma_\chi>0$, $\chi\ne\chi_0$, and with coefficients
\begin{align*}
&\Psi^*_j:\qquad c_j=1-\nu(q,a_j),\qquad r_{j,\chi,\gamma}=2\bar\chi(a_j)\,\frac{u_{\chi,\gamma}}{\rho_\chi},\\
&\pi_j:\qquad c_j=1-\nu(q,a_j),\qquad r_{j,\chi,\gamma}=-2\bar\chi(a_j)\,\frac{1}{\rho_\chi},
\end{align*}
where $u_{\chi,\gamma}=-\bar\rho_\chi/\rho_\chi$ has modulus one. (The common shift $+1$ comes from the bias $-1$ of $\frac{\log x}{\sqrt x}H(x)$ and of $\frac{\log x}{\sqrt x}(\pi(x)-\Li(x))$.) Under LI the non-negative ordinates $\gamma_\chi$, $\chi\ne\chi_0$, are linearly independent over $\mathbb Q$; in particular none of them is $0$. By the vector form of the theorem quoted at the beginning of this section, both vector functions have limiting distributions. These are the laws of
\[
\Bigl(c_j+\operatorname{Re}\sum_{\chi,\gamma}\frac{2\bar\chi(a_j)u_{\chi,\gamma}}{\rho_\chi}Z_{\chi,\gamma}\Bigr)_j\quad\text{and}\quad\Bigl(c_j-\operatorname{Re}\sum_{\chi,\gamma}\frac{2\bar\chi(a_j)}{\rho_\chi}Z_{\chi,\gamma}\Bigr)_j ,
\]
with independent $Z_{\chi,\gamma}$ uniformly distributed on the unit circle. Since $(u_{\chi,\gamma}Z_{\chi,\gamma})$ and $(-Z_{\chi,\gamma})$ are again independent and uniformly distributed, and the rotation $u_{\chi,\gamma}$ does not depend on $j$, the two random vectors have the same distribution $\mu$. For $i\ne j$ some $\chi$ has $\chi(a_i)\ne\chi(a_j)$, so the difference of the $i$-th and $j$-th coordinates has a non-zero coefficient, and its law is absolutely continuous. Hence $\mu(\{\xi_i=\xi_j\})=0$ for $i\ne j$. (The measure $\mu$ itself need not be absolutely continuous: if $r=\varphi(q)$, the random vector lies in the hyperplane $\sum_j\xi_j=0$.)

It remains to pass from the $H(x;q,a_j)$ to the $D(x;q,a_j)$. By Lemma~\ref{lem:HtoD},
\[
D(x;q,a_i)-D(x;q,a_j)=\frac{\sqrt x}{\log x\,\Li(x)}\Bigl(\Psi_i(y)-\Psi_j(y)+O\Bigl(\frac{(\log x)^{6}}{\sqrt x}\Bigr)\Bigr).
\]
Hence for $x$ large, $D(x;q,a_i)>D(x;q,a_j)$ if $\Psi_i(y)-\Psi_j(y)>\lambda$, and $D(x;q,a_i)<D(x;q,a_j)$ if $\Psi_i(y)-\Psi_j(y)<-\lambda$, for any fixed $\lambda>0$. The logarithmic upper and lower densities of the event $\{D(x;q,a_{\tau(1)})>\dots>D(x;q,a_{\tau(r)})\}$ differ from $\mu(\{\xi:\ \xi_{\tau(1)}>\dots>\xi_{\tau(r)}\})$ by at most $\mu(\{\min_{i\ne j}|\xi_i-\xi_j|\le\lambda\})$; here the boundaries $\{\xi_i-\xi_j=\pm\lambda\}$ are $\mu$-null, because the laws of the differences are absolutely continuous. This bound tends to $0$ as $\lambda\to0$, since $\mu(\{\xi_i=\xi_j\})=0$ for $i\ne j$. The same argument, applied to $\pi_j$, gives the density of $\{\pi(x;q,a_{\tau(1)})>\dots>\pi(x;q,a_{\tau(r)})\}$, which is therefore the same number.
\end{proof}

\section{Higher moments: proof of Theorem~\ref{thm:moments}}\label{sec:moments}

Put $F_j(x)=\sum_{p\le x}(\log x/p)^j=\pi(x)M_j(x)$ and $F_j^0(x)=\Li(x)M_j^0(x)$, so that $F_1=F$ and $F_1^0=F_0$. Integration by parts gives
\[
F_j(x)=\int_{2^-}^x\Bigl(\log\frac xt\Bigr)^jd\pi(t)=j\int_2^x\Bigl(\log\frac xt\Bigr)^{j-1}\pi(t)\,\frac{dt}{t},
\]
and the same with $\Li$ in place of $\pi$. Since $dF(t)=\pi(t)\,dt/t$ and $F(2)=0$, a second integration by parts gives, for $j\ge2$,
\begin{equation}\label{eq:Gj}
\mathcal G_j(x):=F_j(x)-F_j^0(x)=j(j-1)\int_2^x\Bigl(\log\frac xt\Bigr)^{j-2}\bigl(F(t)-F_0(t)\bigr)\frac{dt}{t}.
\end{equation}
Moreover, writing $H_j(x)=\pi(x)\bigl(M_j(x)-M_j^0(x)\bigr)$, we have, exactly as in \eqref{eq:H1},
\begin{equation}\label{eq:Hj}
H_j(x)=F_j-M_j^0\pi=\mathcal G_j(x)-M_j^0(x)E(x),
\end{equation}
and $M_j^0(x)=j!+O_j(1/\log x)$, by the expansion of $\Li$.

\begin{lemma}\label{lem:Gj}
Assume the Riemann hypothesis, and let $j\ge1$. Then
\[
\mathcal G_j(x)=-j!\,\frac{\sqrt x}{\log x}\Bigl(2^j+\sum_\rho\frac{x^{i\gamma}}{\rho^{j+1}}\Bigr)+O_j\Bigl(\frac{\sqrt x}{(\log x)^2}\Bigr),
\]
where the series converges absolutely.
\end{lemma}

\begin{proof}
For $j=1$ this is Lemma~\ref{lem:integral}, since $\mathcal G_1=\int_2^xE(t)\,dt/t$. Let $j\ge2$. We insert Lemma~\ref{lem:integral} into \eqref{eq:Gj}. The error term $O(\sqrt t/(\log t)^2)$ contributes $O_j(\sqrt x/(\log x)^2)$. For a complex number $a$ with $\operatorname{Re}a=\frac12$, the substitution $t=xe^{-v}$ gives
\[
j(j-1)\int_2^x\Bigl(\log\frac xt\Bigr)^{j-2}\frac{t^{a-1}}{\log t}\,dt=j(j-1)\,x^a\int_0^{\log(x/2)}\frac{v^{j-2}e^{-av}}{\log x-v}\,dv .
\]
On $v\le\frac12\log x$ we write $(\log x-v)^{-1}=(\log x)^{-1}(1+O(v/\log x))$, and on the remaining range the integrand is $O(v^{j-2}e^{-v/2})$. Since $\int_0^\infty v^{j-2}e^{-av}\,dv=(j-2)!/a^{j-1}$, the expression equals
\[
\frac{j!\,x^a}{a^{j-1}\log x}+O_j\Bigl(\frac{\sqrt x}{(\log x)^2}+x^{1/4}(\log x)^{j}\Bigr),
\]
uniformly in $a$. We apply this with $a=\frac12$, with coefficient $-2$, and with $a=\rho$, with coefficient $-\rho^{-2}$, as given by Lemma~\ref{lem:integral}. The series over $\rho$ may be integrated term by term because $\sum_\rho|\rho|^{-2}<\infty$. Collecting terms, we obtain the lemma.
\end{proof}

\begin{proof}[Proof of Theorem~\ref{thm:moments}]
(a) By Lemma~\ref{lem:Gj}, $\mathcal G_j(x)=O_j(\sqrt x/\log x)$. Hence, by \eqref{eq:Hj}, $H_j(x)=-j!\,E(x)(1+o(1))+O_j(\sqrt x/\log x)$, and Littlewood's theorem gives the assertion as in the proof of Theorem~\ref{thm:main}(b).

(b) By \eqref{eq:Hj}, Lemma~\ref{lem:Gj} and the argument of Lemma~\ref{lem:normal},
\[
\frac{\log x}{j!\sqrt x}H_j(x)=-(2^j-1)+2\operatorname{Re}\sum_{0<\gamma\le T}x^{i\gamma}\,w_{j,\gamma}+\varepsilon_{j,T}(y),\qquad w_{j,\gamma}=\frac{1-\rho^{-j}}{\rho},
\]
where the mean square of $\varepsilon_{j,T}$ is $\ll_j(\log T)^2/T$. Here $-M_j^0E$ contributes $1+2\operatorname{Re}\sum x^{i\gamma}/\rho$ up to such errors, and $\mathcal G_j$ contributes $-2^j-2\operatorname{Re}\sum x^{i\gamma}/\rho^{j+1}$. By the theorem quoted at the beginning of Section~\ref{sec:mirror}, a limiting distribution exists. Its mean is the constant term $-(2^j-1)$, since the oscillating terms have mean zero and the bounded mean square gives uniform integrability. The limiting distribution of $\frac{\log x}{\sqrt x}(\pi(x)-\li(x))$ has mean $-1$ \cite{RS}, so the two distributions differ when $2^j-1\ne1$, that is, when $j\ge2$. For $j=1$ they coincide under LI by the proof of Theorem~\ref{thm:skewes}.

(c) Under LI the limiting distribution is that of $-(2^j-1)+X_j$, where
\[
X_j=\sum_{\gamma>0}r_{j,\gamma}\cos\vartheta_\gamma,\qquad r_{j,\gamma}=2|w_{j,\gamma}|,
\]
with independent $\vartheta_\gamma$ uniformly distributed on $[0,2\pi)$. This law is absolutely continuous, since each summand has the arcsine law, so $\delta_j=\mathbb P(X_j>2^j-1)$ exists; it is positive by the lower bound proved below. For $j=1$, Theorem~\ref{thm:skewes} identifies $\delta_1$. Since $|w_{j,\gamma}-1/\rho|=|\rho|^{-j-1}\le|\rho|^{-2}$ and $|\rho|=\gamma+O(\gamma^{-1})$, we have
\[
r_{j,\gamma}=\frac2\gamma+O\Bigl(\frac1{\gamma^2}\Bigr)
\]
uniformly in $j$ and $\gamma$. We now estimate $\mathbb P(X_j>V)$ for large $V$, following the method of Montgomery \cite{Montgomery}. The Riemann--von Mangoldt formula $N(T)=\frac{T}{2\pi}\log\frac{T}{2\pi e}+O(\log T)$ gives, by partial summation, for $U\ge U_0$,
\[
N(U)\le U\log U,\qquad \sum_{0<\gamma\le U}\frac1\gamma=\frac{(\log U)^2}{4\pi}+O(\log U),\qquad \sum_{\gamma>U}\frac1{\gamma^2}\ll\frac{\log U}{U}.
\]
Let $0<\varepsilon<\frac14$.

\emph{Upper bound.} For $\lambda>0$ we have $\mathbb E\,e^{\lambda r\cos\vartheta}=I_0(\lambda r)$, where $I_0(t)=\sum_{n\ge0}(t/2)^{2n}/(n!)^2$. We have $I_0(t)\le e^{t}$ and $I_0(t)\le\sum_n(t^2/4)^n/n!=e^{t^2/4}$ for $t\ge0$. By independence, splitting at $\gamma=\lambda$,
\[
\log\mathbb E\,e^{\lambda X_j}\le\sum_{\gamma\le\lambda}\lambda r_{j,\gamma}+\sum_{\gamma>\lambda}\frac{\lambda^2r_{j,\gamma}^2}{4}\le\frac{\lambda(\log\lambda)^2}{2\pi}+O(\lambda\log\lambda),
\]
uniformly in $j$. Hence $\mathbb P(X_j>V)\le\exp\bigl(\frac{\lambda(\log\lambda)^2}{2\pi}(1+\eta)-\lambda V\bigr)$, where $\eta=\eta(\lambda)\to0$ as $\lambda\to\infty$, uniformly in $j$. Choose $\lambda=\exp\bigl((1-\varepsilon)\sqrt{2\pi V}\bigr)$. Then $\frac{(\log\lambda)^2}{2\pi}(1+\eta)=(1-\varepsilon)^2(1+\eta)V\le(1-\varepsilon)V$ once $V$ is large, and therefore $\mathbb P(X_j>V)\le\exp(-\varepsilon\lambda V)$, that is,
\[
\log\log\frac{1}{\mathbb P(X_j>V)}\ge(1-\varepsilon)\sqrt{2\pi V}.
\]

\emph{Lower bound.} Let $U\ge U_0$ and split $X_j=S+Y$, where $S$ is the sum over $\gamma\le U$ and $Y$ the sum over $\gamma>U$. These are independent. With probability $(\varepsilon/\pi)^{N(U)}\ge\exp(-U\log U\log(\pi/\varepsilon))$ every $\vartheta_\gamma$ with $\gamma\le U$ lies in $[-\varepsilon,\varepsilon]$, and then $\cos\vartheta_\gamma\ge1-\varepsilon$ and
\[
S\ge(1-\varepsilon)\sum_{\gamma\le U}r_{j,\gamma}\ge(1-\varepsilon)^2\,\frac{(\log U)^2}{2\pi}
\]
for $U$ large, uniformly in $j$. On the other hand $\mathbb EY=0$ and $\mathbb EY^2=\frac12\sum_{\gamma>U}r_{j,\gamma}^2\ll(\log U)/U$, so $\mathbb P(Y>-1)\ge\frac12$ for $U$ large, by Chebyshev's inequality. Hence
\[
\mathbb P\Bigl(X_j>(1-\varepsilon)^2\frac{(\log U)^2}{2\pi}-1\Bigr)\ge\frac12\exp\Bigl(-U\log U\log\frac\pi\varepsilon\Bigr).
\]
Given a large $V$, choose $U$ with $\log U=(1-\varepsilon)^{-1}\sqrt{2\pi(V+1)}$, so that the left-hand argument above equals $V$. Then
\[
\log\log\frac{1}{\mathbb P(X_j>V)}\le\log U+\log\log U+O_\varepsilon(1)\le(1+2\varepsilon)\sqrt{2\pi V}.
\]

Both bounds hold uniformly in $j\ge1$ for all sufficiently large $V$. Taking $V=2^j-1$, so that $\delta_j=\mathbb P(X_j>V)$, and letting $\varepsilon\to0$, we obtain $\log\log(1/\delta_j)=(1+o(1))\sqrt{2\pi(2^j-1)}=(\sqrt{2\pi}+o(1))2^{j/2}$. For $j=1$ the same argument gives a weak form of Montgomery's estimate \cite{Montgomery}, which in fact gives $\log(1/\mathbb P(X_1>V))\asymp\sqrt V\exp(\sqrt{2\pi V})$; the asymptotic form of the large deviations of this distribution is due to Monach \cite{Monach}, and races between several residue classes were treated by Lamzouri \cite{Lamzouri}.
\end{proof}

\section{Function fields: proof of Theorem~\ref{thm:FF}}\label{sec:FF}

We follow the notation of Cha \cite{Cha} where possible. Let $\Phi=\Phi(M)$ be the number of residues prime to $M$, and for $N\ge1$ let $\pi^*_q(N)$ and $\pi^*_q(N;M,A)$ be the numbers of monic irreducible $P$ of degree exactly $N$, and of those with $P\equiv A\pmod M$. Put
\[
\Pi(X)=\sum_{N\le X}\pi^*_q(N),\qquad \Pi_A(X)=\pi_q(X;M,A)=\sum_{N\le X}\pi^*_q(N;M,A).
\]
For a non-principal character $\chi$ modulo $M$ write $L(u,\chi)=\prod_\nu(1-\gamma_\nu u)$. By Weil's theorem each inverse zero has $|\gamma_\nu|=1$ or $|\gamma_\nu|=\sqrt q$, and we write $\gamma=\sqrt q\,e^{i\alpha_\gamma}$ in the latter case. Let $c(M,A)=-1+\#\{B:\ B^2\equiv A\pmod M\}$, and let $B_q(X)=\sqrt q/(q-1)$ for odd $X$ and $B_q(X)=q/(q-1)$ for even $X$. Cha's explicit formula \cite[Theorem~2.5]{Cha} states that
\begin{equation}\label{eq:Cha}
\begin{split}
\frac{X}{q^{X/2}}\bigl(\Phi\,\Pi_A(X)-\Pi(X)\bigr)&=-c(M,A)B_q(X)\\
&\qquad-\sum_{\chi\ne\chi_0}\bar\chi(A)\sum_{|\gamma|=\sqrt q}e^{i\alpha_\gamma X}\frac{\gamma}{\gamma-1}+o(1),
\end{split}
\end{equation}
where $\gamma$ runs over the inverse zeros of modulus $\sqrt q$ of $L(u,\chi)$. These are finite sums.

The moment of the system about the point $X$ is
\[
F_A(X)=\sum_{N\le X}(X-N)\,\pi^*_q(N;M,A)=\sum_{Y=1}^{X-1}\Pi_A(Y),
\]
the analogue of the identity $F(x)=\int_2^x\pi(t)\,dt/t$, and $D_q(X;M,A)=F_A(X)/\Pi_A(X)$. Define $F(X)=\sum_{Y<X}\Pi(Y)$ and $D_q(X)=F(X)/\Pi(X)$ in the same way for all monic irreducibles. Since $\pi^*_q(N)=q^N/N+O(q^{N/2}/N)$, we have $\Pi(X)=\frac{q}{q-1}\frac{q^X}{X}(1+O(\frac1X))$ and
\[
D_q(X)=\frac{1}{q-1}+O\Bigl(\frac1X\Bigr).
\]
This is the analogue of $D(x)\to1$: seen from the point $X$, the masses are distributed according to the geometric law with ratio $1/q$, whose mean is $1/(q-1)$.

\begin{lemma}\label{lem:FF}
Put $H_A(X)=\Pi_A(X)\bigl(D_q(X;M,A)-D_q(X)\bigr)$. Then
\[
\begin{aligned}
\frac{q-1}{\sqrt q}\cdot\frac{X}{q^{X/2}}\,\Phi\,H_A(X)&=-c(M,A)B_q(X+1)\\
&\qquad-\sum_{\chi\ne\chi_0}\bar\chi(A)\sum_{|\gamma|=\sqrt q}e^{i\alpha_\gamma X}\frac{\gamma}{\gamma-1}\,v_\gamma+o(1),
\end{aligned}
\]
where
\[
v_\gamma=\frac{q-1}{\sqrt q}\Bigl(\frac{1}{\gamma-1}-\frac{1}{q-1}\Bigr)=\frac{q-\gamma}{\sqrt q\,(\gamma-1)}\qquad\text{satisfies}\qquad|v_\gamma|=1 .
\]
\end{lemma}

\begin{proof}
Since $F=D_q\Pi$, we have $\Phi H_A=(\Phi F_A-F)-D_q(X)(\Phi\Pi_A-\Pi)$. Write $R(Y)$ for the left-hand side of \eqref{eq:Cha}. Then
\[
\Phi F_A(X)-F(X)=\sum_{Y<X}\frac{q^{Y/2}}{Y}\,R(Y).
\]
For $|\gamma|=\sqrt q$ we have $\sum_{Y<X}\gamma^Y/Y=\frac{\gamma^X}{X}\cdot\frac{1}{\gamma-1}+o(q^{X/2}/X)$, by summing the geometric series backwards from $Y=X-1$. Similarly, splitting according to the parity of $Y$,
\[
\begin{aligned}
\frac{X}{q^{X/2}}\sum_{Y<X}\frac{q^{Y/2}}{Y}B_q(Y)&=\sum_{j\ge1}q^{-j/2}B_q(X-j)+o(1)\\
&=\frac{\sqrt q}{q-1}B_q(X+1)+\frac{B_q(X)}{q-1}+o(1),
\end{aligned}
\]
the right-hand side being $2q/(q-1)^2$ for even $X$ and $\sqrt q(q+1)/(q-1)^2$ for odd $X$. The terms $o(1)$ in \eqref{eq:Cha} contribute $o(q^{X/2}/X)$. Subtracting $D_q(X)(\Phi\Pi_A-\Pi)$, with $D_q(X)=\frac1{q-1}+O(\frac1X)$, the terms $B_q(X)/(q-1)$ cancel, and the zero $\gamma$ acquires the factor $\frac1{\gamma-1}-\frac1{q-1}$. Finally, with $\gamma=\sqrt q\,e^{i\alpha_\gamma}$,
\[
|q-\gamma|^2=q^2-2q\sqrt q\cos\alpha_\gamma+q=q\,|\gamma-1|^2,
\]
so $|v_\gamma|=1$.
\end{proof}

\begin{proof}[Proof of Theorem~\ref{thm:FF}]
Let $\mathbf h(X)$ be the vector of the left-hand sides of Lemma~\ref{lem:FF} for $A=A_1,\dots,A_r$, and let $\mathbf e(X)$ be the vector of the left-hand sides of \eqref{eq:Cha}. The conjugates of the inverse zeros of $L(u,\chi)$ are those of $L(u,\bar\chi)$, so the $\gamma$ occur in conjugate pairs; under GSH none of them is real. Moreover, by GSH the numbers $\alpha_\gamma$ with $0<\alpha_\gamma<\pi$, together with $2\pi$, are linearly independent over $\mathbb Q$. By Kronecker's theorem, as $X$ runs over the integers, the point $(X\alpha_\gamma/2\pi)_\gamma$ of the torus together with the parity of $X$ is equidistributed with respect to the product of the Haar measure and the uniform measure on $\{0,1\}$. Hence, for each parity $\kappa$, the vectors $\mathbf e(X)$ and $\mathbf h(X)$, with $X$ restricted to the integers of parity $\kappa$, have limiting distributions. These are the laws of
\[
\begin{gathered}
\Bigl(-c(M,A_j)B_q(\kappa)-\operatorname{Re}\sum_\gamma W_{j,\gamma}\Bigr)_j\\ \text{and}\\ \Bigl(-c(M,A_j)B_q(\kappa+1)-\operatorname{Re}\sum_\gamma v_\gamma W_{j,\gamma}\Bigr)_j ,
\end{gathered}
\]
where $\gamma$ runs over the inverse zeros with $0<\alpha_\gamma<\pi$ of all $L(u,\chi)$, and $W_{j,\gamma}=2\bar\chi(A_j)\frac{\gamma}{\gamma-1}e^{i\vartheta_\gamma}$, $\chi$ being the character with the inverse zero $\gamma$, with independent $\vartheta_\gamma$ uniformly distributed on $[0,2\pi)$. The pair of terms for $\gamma$ and $\bar\gamma$ is combined into a real part, because $v_{\bar\gamma}=\overline{v_\gamma}$. Since $|v_\gamma|=1$ and $v_\gamma$ does not depend on $j$, the variables $e^{i\vartheta_\gamma}v_\gamma$ are again independent and uniform. Therefore the limiting law of $\mathbf h$ over the integers of parity $\kappa$ equals the limiting law of $\mathbf e$ over the integers of parity $\kappa+1$.

By the hypothesis of the theorem, for $i\ne j$ the difference of the $i$-th and $j$-th coordinates of either random vector contains a term $\operatorname{Re}(cZ)$ with $c\ne0$ and $Z$ uniform on the unit circle, independent of the other terms. Its law is therefore continuous, and the hyperplanes $\{x_i=x_j\}$ are null for each of the four laws. The existence of the limiting laws over each parity follows from the Kronecker argument above (compare \cite[Theorem~3.2]{Cha}); GSH is \cite[Definition~3.3]{Cha}. Finally, $\Phi\Pi_A(X)=\Pi(X)(1+O(Xq^{-X/2}))$ and $H_A(X)\ll q^{X/2}/X$, so
\[
D_q(X;M,A_i)-D_q(X;M,A_j)=\frac{\Phi}{\Pi(X)}\bigl(H_{A_i}(X)-H_{A_j}(X)\bigr)+O\Bigl(\frac{1}{\Pi(X)}\Bigr),
\]
so the ordering of the $D_q(X;M,A_j)$ agrees with that of the coordinates of $\mathbf h(X)$ whenever these coordinates differ pairwise by more than a fixed $\lambda>0$ and $X$ is large. As in the proof of Theorem~\ref{thm:race}, the set of $X$ where some two coordinates differ by at most $\lambda$ has upper density tending to $0$ with $\lambda$, since the hyperplanes $\{x_i=x_j\}$ are null. This gives both assertions of the theorem.
\end{proof}

\section*{Concluding remark: other weights}

The centre of mass corresponds to the weight $w(v)=v$ in the sums $\sum_{p\le x}w(\log(x/p))$. For a smooth weight $w$ the explicit formula attaches to the zero $\rho$ a weight given by the Mellin transform of $w$, and to the squares of primes a corresponding bias; the weights $w(v)=v^j$ of Theorem~\ref{thm:moments} show how both change with $w$. The general frameworks of Akbary, Ng and Shahabi \cite{ANS}, Devin \cite{Devin} and Fiorilli and Martin \cite{FM}, and the work of Meng \cite{Meng} on products of $k$ primes, apply to such sums. We do not know for which weights the limiting distribution coincides with that of $\pi(x)-\li(x)$, nor for which a pointwise statement such as Theorem~\ref{thm:mirror} holds.

\end{document}